\documentclass[12pt]{article}
\usepackage{amssymb}
\usepackage{amsthm}
\usepackage{bbold}
\input diagxy

\def\spaces{\mathbb{T}{\rm op}}
\def\topos{\mathbb{T}{\rm opos}}
\def\loc{\mathbb{L}{\rm oc}}
\def\cat{\mathbb{C}{\rm at}}
\def\pos{\mathbb{P}{\rm os}}
\def\ring{\mathbb{R}{\rm ing}}
\def\quant{\mathbb{Q}{\rm uant}}
\def\span{\mathbb{S}{\rm pan}}
\def\cospan{\mathbb{C}{\rm ospan}}
\def\tov{\hbox{$\to<150>$\hskip -.2in \raisebox{1.7pt}{\tiny$\bullet$} \hskip .1in}}
\def\id{{\rm id}}
\def\idv{\id^{\hbox{\tiny$\bullet$}}}
\def\dotv{\raisebox{1.5pt}{\tiny$\bullet$}}

\title{Span, Cospan, and Other Double Categories}
\vfill 
\author{Susan Niefield}

\begin{document}

\maketitle

\newtheorem{lem}{Lemma}[section]
\newtheorem{prop}[lem]{Proposition}
\newtheorem{cor}[lem]{Corollary}
\newtheorem{thm}[lem]{Theorem}

\theoremstyle{definition}
\newtheorem{ex}[lem]{Example}

\begin{abstract}
Given a double category $\mathbb D$ such that $\mathbb D_0$ has pushouts, we characterize oplax/lax adjunctions $\mathbb D\two/->`<-/<150>\cospan(\mathbb D_0)$ for which the right adjoint is normal and restricts to the identity on $\mathbb D_0$, where $\cospan(\mathbb D_0)$ is the double category on $\mathbb D_0$ whose vertical morphisms are cospans.  We show that such a pair exists if and only if $\mathbb D$ has companions, conjoints, and 1-cotabulators. The right adjoints are induced by the companions and conjoints, and the left adjoints by the 1-cotabulators.  The notion of a 1-cotabulator is a common generalization of the symmetric algebra of a module and Artin-Wraith glueing of toposes, locales, and topological spaces.
\end{abstract}

\section{Introduction}
Double categories, first introduced by Ehresmann \cite{Ehres}, provide a setting in which one can simultaneously consider two kinds of morphisms (called {\it horizontal} and {\it vertical} morphisms).  Examples abound in many areas of mathematics.  There are double categories whose objects are sets, rings, categories, posets, topological spaces, locales, toposes, quantales, and more.

\smallskip
There are general examples, as well.  If $\cal D$ is a category with pullbacks, then there is a double category $\span({\cal D})$ whose objects and horizontal morphisms are those of $\cal D$, and vertical morphisms $X_0\tov X_1$ are spans, i.e., morphisms  $X_0\toleft<125>X\to<125>X_1$ of $\cal D$, with vertical compositions via pullback.  If $\cal D$ is a category with pushouts, then $\cospan(\cal D)$ is defined dually, in the sense that vertical morphisms $X_0\tov X_1$ are cospans $X_0\to<125>X\toleft<125>X_1$ with vertical compositions via pushout.  Moreover, if $\cal D$ has both, then pushout of spans and pullback of cospans induce an oplax/lax adjunction (in the sense of Par\'{e} \cite{Yoneda})  $$\span({\cal D})\two/->`<-/^F_G\cospan({\cal D})$$ which restricts to the identity on the horizontal category $\cal D$.

\smallskip
Now, suppose $\cal D$ is a category with pushouts, and we replace $\span({\cal D})$ by a double category $\mathbb D$ whose horizontal category is also $\cal D$.  Then several questions arise.  Under what conditions on $\mathbb D$ is there an oplax/lax adjunction $\mathbb D\two/->`<-/<150>\cospan({\cal D})$ which restricts to the identity on  $\cal D$?  In particular, if $\mathbb D$ has cotabulators, then there is an induced oplax functor $F\colon \mathbb D \to<125> \cospan({\cal D})$, and so one can ask when this functor $F$ has a right adjoint.  Similarly,  if $\mathbb D$ has companions and  conjoints, there is a normal lax functor $G\colon \cospan({\cal D})\to<125> \mathbb D$ which takes a cospan $X_0\to<125>^{c_0}X\toleft<125>^{c_1}X_1$ to the composite $X_0{\hbox{$\to<150>^{{c_0}_*}$\hskip -.2in \raisebox{1.7pt}{\tiny$\bullet$} \hskip .1in}}X{\hbox{$\to<150>^{c_1^*}$\hskip -.2in \raisebox{1.7pt}{\tiny$\bullet$} \hskip .1in}}X_1$, and so one can ask when the induced functor $G$ has a left adjoint.  

\smallskip
We will see that these questions are related.  In particular, there is an oplax/lax adjunction $$\mathbb D\two/->`<-/^F_G\cospan({\cal D})$$ such that $G$ is normal and restricts to the identity on $\cal D$ precisely when $\mathbb D$ has 1-cotabulators, conjoints, and companions  (where for 1-cotabulators we drop the tetrahedron condition in the definition of cotabulator).  Moreover, if $\cal D$ has pullbacks, then the result dualizes to show that there is oplax/lax adjunctions $\span({\cal D}) \two/->`<-/<150>\mathbb D$ whose left adjoint is opnormal and restricts to the identity on $\cal D$ precisely when $\mathbb D$ has 1-tabulators, conjoints, and companions. 

\smallskip
The double categories mentioned above all have companions, conjoints, and 1-cotabulators, and the functors $F$ and $G$ are related to familiar constructions.  In the double category of commutative rings (as well as, quantales), the functor $F$ is given by the symmetric algebra algebra on a bimodule  and $G$ is given by restriction of scalars.  For categories (and posets), $F$ is the collage construction.   In the case of topological spaces, locales, and toposes, the functor $F$ uses Artin-Wraith glueing.

\smallskip
The paper proceeds as follows.  We begin in Section~2 with the  double categories under consideration, followed by a review of companions and conjoints in Section~3.  The notion of 1-tabulators (duallly, 1-cotabulators) is then introduced in Section~4.  After a brief discussion of oplax/lax adjunctions in Section~5, we present our characterization (Theorem~\ref{adj}) of those of the form $\mathbb D\two/->`<-/<150>\cospan({\cal D})$ such that right adjoint is normal and restricts to the identity on $\cal D$.  Along the way, we obtain a possibly new characterization (Proposition~\ref{framed}) of double categories with companions and conjoints, in the case where the horizontal category $\cal D$ has pushouts, as those for which the identity functor on $\cal D$ extends to a normal lax functor ${\cospan}(\cal D)\to<125> \mathbb D$.  We conclude with the dual (Corollary~\ref{adjop}) classification of oplax/lax adjunctions $\span({\cal D}) \two/->`<-/<150>\mathbb D$  left adjoint is opnormal and restricts to the identity on $\cal D$.

\section{The Examples of Double Categories}

Following Par\'{e} \cite{GP99,Yoneda} and Shulman \cite{Shulman}, we define a {\it double category} $\mathbb D$ to be a weak internal category 
$$\bfig 
\morphism(0,0)<450,0>[\mathbb D_1\times_{\mathbb D_0}\mathbb D_1`\mathbb D_1;c]
\morphism(450,0)|a|/@{>}@<5pt>/<350,0>[\mathbb D_1`\mathbb D_0;d_0]
\morphism(800,0)|m|<-350,0>[\mathbb D_0`\mathbb D_1;\Delta]
\morphism(450,0)|b|/@{>}@<-5pt>/<350,0>[\mathbb D_1`\mathbb D_0;d_1]
\efig$$ 
in $\rm CAT$.   It consists of objects (those of $\mathbb D_0$), two types of morphisms: horizontal (those of $\mathbb D_0$) and vertical (objects of $\mathbb D_1$ with domain and codomain given by $d_0$ and $d_1$), and cells (morphisms of $\mathbb D_1$) denoted by
$$\bfig
\square<350,350>[X_0`Y_0`X_1`Y_1;f_0`m`n`f_1]
\place(0,200)[\hbox{\tiny$\bullet$}]
\place(350,200)[\hbox{\tiny$\bullet$}]
\place(175,175)[\varphi]
\efig\eqno{(\star)}$$
Composition and identity morphisms are given horizontally in $\mathbb D_0$  and vertically via $c$ and $\Delta$, respectively.  

\smallskip
The objects, horizontal morphisms, and  {\it special} cells (i.e., ones in which the vertical morphisms are identities) form a 2-category called the {\it horizontal 2-category} of $\mathbb D$.  Since $\mathbb D$ is a weak internal category in $\rm CAT$, the associativity and identity axioms for vertical morphisms hold merely up to coherent isomorphism, and so we get an  analogous {\it vertical bicategory}.  When these isomorphisms are identities, we say that $\mathbb D$ is a {\it strict} double category.

\medskip
The following double categories are of interest in this paper.

\begin{ex} $\spaces$ has topological spaces as objects and continuous maps as horizontal morphisms.  Vertical morphisms $X_0\tov X_1$ are finite intersection-preserving maps ${\cal O}(X_0)\tov {\cal O}(X_1)$ on the open set lattices, and there is a cell of the form $(\star)$ if and only if $f_1^{-1}n\subseteq mf_0^{-1}$.
\end{ex}

\begin{ex} $\loc$ has locales as objects, locale morphisms (in the sense of \cite{SS}) as horizontal morphisms, and
 finite meet-preserving maps as vertical morphisms. There is a cell of the form $(\star)$ if and only if $f_1^*n\le mf_0^*$.
\end{ex}

\begin{ex} $\topos$ has Grothendieck toposes as objects, geometric morphisms (in the sense of \cite{TT}) as horizontal morphisms, and natural transformations $f_1^*n\to<125>mf_0^*$ as cells of the form  $(\star)$.
\end{ex}

\begin{ex} $\cat$ has small categories as objects and functors as horizontal morphisms. Vertical morphisms $m\colon X_0\tov X_1$ are profunctors (also known as distributors and relators), i.e., functors  $m\colon X_0^{op}\times X_1\to<125>{\rm Sets}$, and  natural transformations $m\to<125>n(f_0-,f_1-)$ are cells of the form $(\star)$.
\end{ex}

\begin{ex} $\pos$ has partially-ordered sets as objects and order-preserving maps as horizontal morphisms.  Vertical morphisms $m\colon X_0\tov X_1$ are order ideals $m\subseteq  X_0^{op}\times X_1$, and there is a cell of the form $(\star)$ if and only if $(x_0,x_1)\in m \Rightarrow (f_0(x_0),f_1(x_1))\in n$.
\end{ex}

\begin{ex} For a category $\cal D$ with pullbacks, the span double category  $\span({\cal D})$ has objects and horizontal morphisms in $\cal D$, and vertical morphisms which are spans in  $\cal D$, with composition defined via pullback and the identities $\idv\colon X\tov X$ given by $X\to<150>^{\id_X} X\to/<-/<150>^{\id_X} X$. The cells $m\to<125>n$  are commutative diagrams in $\cal D$ of the form
$$\bfig
\Ctriangle/<-``->/<350,200>[X_0`X`X_1;m_0``m_1]
\morphism(0,200)|m|<600,0>[X`Y;f]
\morphism(350,400)|a|<600,0>[X_0`Y_0;f_0]
\morphism(350,0)|b|<600,0>[X_1`Y_1;f_1]
\Ctriangle(600,0)|rrr|/<-``->/<350,200>[Y_0`Y`Y_1;n_0``n_1]
\efig$$
In particular, $\span({\rm Set})$ is the double category $\mathbb S{\rm et}$ considered by Par\'{e} in \cite{Yoneda}, see also \cite{DPPSpan}.
\end{ex}

\begin{ex} $\cospan({\cal D})$ is defined dually, for a category $\cal D$ with pushouts.  In particular, $\span({\rm Top})$ is the double category used by Grandis \cite{cospan1,cospan2} in his study of 2-dimensional topological quantum field theory.
\end{ex}

\begin{ex} For a symmetric monoidal category $\cal V$ with coequalizers, the double category $\mathbb M{\rm od}({\cal V})$ has commutative monoids in $\cal V$ as objects and monoid homomorphisms as horizontal morphisms.  Vertical morphisms from $X_0$ to $X_1$  are $(X_0,X_1)$-bimodules, with composition via tensor product, and cells are bimodule homomorphisms.  Special cases include the double category $\ring$ of commutative rings with identity and the double category $\quant$ of commutative unital quantales.
\end{ex}

\section{Companions and Conjoints}

Recall \cite{GP04} that companions and conjoints in a double category are defined as follows.  Suppose $f\colon X\to<125>Y$ is a  horizontal morphism.  
A {\it companion} for $f$ is a vertical morphism $f_*\colon X\tov Y$ together with cells 
$$\bfig
\square(0,0)<300,300>[X`X`X`Y;\id_X`\idv_X`f_*`f]
\place(0,180)[\hbox{\tiny${\bullet\atop{ \ }}$}]
\place(300,180)[\hbox{\tiny${\bullet\atop{ \ }}$}]
\place(150,155)[\eta]
\square(1000,0)<300,300>[X`Y`Y`Y;f`f_*`\idv_Y`\id_Y]
\place(1000,180)[\hbox{\tiny${\bullet\atop{ \ }}$}]
\place(1300,180)[\hbox{\tiny${\bullet\atop{ \ }}$}]
\place(1150,155)[\varepsilon]
\efig$$
whose horizontal and vertical compositions are identity cells.   A {\it conjoint} for $f$ is a vertical morphism  morphism $f^*\colon Y\tov X$ together with cells 
$$\bfig
\square(0,0)<300,300>[X`Y`X`X;f`\id_X`f^*`\id_X]
\place(0,180)[\hbox{\tiny${\bullet\atop{ \ }}$}]
\place(300,180)[\hbox{\tiny${\bullet\atop{ \ }}$}]
\place(150,155)[\alpha]
\square(1000,0)<300,300>[Y`Y`X`Y;\id_Y`f^*`\idv_Y`f]
\place(1000,180)[\hbox{\tiny${\bullet\atop{ \ }}$}]
\place(1300,180)[\hbox{\tiny${\bullet\atop{ \ }}$}]
\place(1150,155)[\beta]
\efig$$
whose horizontal and vertical compositions are identity cells.  We say {\it $\mathbb D$ has companions and conjoints} if every horizontal morphism has a companion and a conjoint.  Such a double category is also known as a {\it framed bicategory} \cite{Shulman}.

\smallskip
If $f$ has a companion $f_*$, then one can show that there is a bijection between cells of the following form 
$$\bfig
\square(0,150)/`->`->`->/<600,300>[\cdot`Y`\cdot`\cdot;`m`n`h]
\morphism(0,450)<300,0>[\cdot`X;g]
\morphism(300,450)<300,0>[X`Y;f]
\place(0,300)[\hbox{\tiny${\bullet\atop{ \ }}$}]
\place(600,300)[\hbox{\tiny${\bullet\atop{ \ }}$}]
\place(300,300)[\varphi]
\square(1200,0)/->`->``->/<350,600>[\cdot`X`\cdot`\cdot;g`m``h]
\morphism(1550,600)|r|<0,-300>[X`Y;f_*]
\morphism(1550,300)|r|<0,-300>[Y`\cdot;n]
\place(1200,300)[\hbox{\tiny${\bullet\atop{ \ }}$}]
\place(1550,150)[\hbox{\tiny${\bullet\atop{ \ }}$}]
\place(1550,450)[\hbox{\tiny${\bullet\atop{ \ }}$}]
\place(1375,300)[\psi]
\efig$$
Similarly, if $f$ has a conjoint $f^*$, then there is a bijection between cells
$$\bfig
\square(0,150)/->`->`->`/<600,300>[\cdot`\cdot`\cdot`Y;g`m`n`]
\morphism(0,150)|b|<300,0>[\cdot`X;h]
\morphism(300,150)|b|<300,0>[X`Y;f]
\place(0,300)[\hbox{\tiny${\bullet\atop{ \ }}$}]
\place(600,300)[\hbox{\tiny${\bullet\atop{ \ }}$}]
\place(300,320)[\varphi]
\square(1200,0)/->`->``->/<350,600>[\cdot`\cdot`\cdot`X;g`m``h]
\morphism(1550,600)|r|<0,-300>[\cdot`Y;n]
\morphism(1550,300)|r|<0,-300>[Y`X;f^*]
\place(1200,300)[\hbox{\tiny${\bullet\atop{ \ }}$}]
\place(1550,150)[\hbox{\tiny${\bullet\atop{ \ }}$}]
\place(1550,450)[\hbox{\tiny${\bullet\atop{ \ }}$}]
\place(1375,300)[\psi]
\efig$$
There are two other cases of this process (called {\it vertical flipping} in \cite{GP04}) which we do not recall here as they will not be used in the following. 

\smallskip
All of the double categories mention in the previous section have well known companions and conjoints.  In $\spaces$, $\loc$, and $\topos$, the companion and conjoint of $f$ are the usual maps denoted by $f_*$ and $f^*$.  For $\cat$, they are the profunctors defined by $f_*(x,y)=Y(fx,y)$ and $f^*(y,x)=Y(y,fx)$, and analogously, for $\pos$.   If $\cal V$ is a symmetric monoidal category and $f\colon X\to<125> Y$  is a monoid homomorphism, then $Y$ becomes an $(X,Y)$-bimodule and a $(Y,X)$-bimodule via $f$, and so $Y$ is both a companion and conjoint for $f$.  Finally, for $\span({\cal D})$  (respectively, $\cospan({\cal D})$) the companion and conjoint of $f$ are the span (respectively, cospan) with $f$ as one leg and the appropriate identity morphism as the other.

\section{1-Tabulators and 1-Cotabulators}

Tabulators in double categories were defined as follows in \cite{GP99} (see also \cite{Yoneda}).  Suppose $\mathbb D$ is a double category and $m\colon X_0\tov X_1$ is a vertical morphism in $\mathbb D$.  A {\it tabulator} of $m$ is an object $T$ together with a cell 
$$\bfig
\Ctriangle/<-`->`->/<350,200>[X_0`T`X_1;`m`]
\place(220,200)[\tau]
\place(350,200)[\hbox{\tiny${\bullet\atop{ \ }}$}]
\efig$$
such that for any cell 
$$\bfig
\Ctriangle/<-`->`->/<350,200>[X_0`Y`X_1;`m`]
\place(220,200)[\varphi]
\place(350,200)[\hbox{\tiny${\bullet\atop{ \ }}$}]
\efig$$
there exists a unique morphism $f\colon Y\to<125>T$ such that $\tau f=\varphi$, and for any commutative tetrahedron of cells
$$\bfig
\square<350,350>[Y_0`X_0`Y_1`X_1;`n`m`]
\place(0,200)[\hbox{\tiny$\bullet$}]
\place(350,200)[\hbox{\tiny$\bullet$}]
\morphism(0,350)<350,-350>[Y_0`X_1;]
\morphism(0,0)<350,350>[Y_1`X_0;]
\efig$$
there is a unique cell $\xi$ such that 
$$\bfig
\Dtriangle/->`->`<-/<300,200>[Y_0`T`Y_1;n``]
\Ctriangle(300,0)/<-`->`->/<300,200>[X_0`T`X_1;`m`]
\place(500,200)[\tau]
\place(100,200)[\xi]
\place(0,200)[\hbox{\tiny${\bullet\atop{ \ }}$}]
\place(600,200)[\hbox{\tiny${\bullet\atop{ \ }}$}]
\efig$$
gives the tetrahedron in the obvious way. 

\smallskip
Tabulators (and their duals cotabulators) arise in the next section, but we do not use the tetrahedron property in any of our proofs or constructions.  Thus, we drop this condition in favor of a weaker notion which we call a {\it 1-tabulator} (and dually, {\it 1-cotabulator}) of a vertical morphism.  When the tetrahedron condition holds, we call these tabulators (respectively, cotabulators) {\it strong}.  

\smallskip
It is easy to show that the following proposition gives an alternative definition in terms of adjoint functors.

\begin{prop}\label{tab} A double category $\mathbb D$ has 1-tabulators (respectively, 1-cotabulators) if and only if $\Delta\colon \mathbb D_0\to<125> \mathbb D_1$ has a right (respectively, left) adjoint which we denote by  $\Sigma$ (respectively, $\Gamma$).
\end{prop}

\begin{cor}\label{tabcor} $\mathbb M{\rm od}({\cal V})$ does not have 1-tabulators.
\end{cor} 

\noindent{\bf Proof}.  Suppose $({\cal V}, \otimes, I)$ is a symmetric monoidal category.  Then $I$ is an initial object $\mathbb M{\rm od}({\cal V})_0$, which is the category of commutative monoids in $\cal V$.  Since $\Delta I$ is not an initial object of $\mathbb M{\rm od}({\cal V})_1$, we know $\Delta$ does not have a right adjoint, and so the result follows from Proposition~\ref{tab}.
\hfill$\Box$

\medskip
The eight examples under consideration have 1-cotabulators, and we know that 1-tabulators exist in all but one (namely,  $\mathbb M{\rm od}({\cal V})$).  In fact, we will prove a general proposition that gives the existence of 1-tabulators from a property of 1-cotabulators (called {\it $\mathbb2$-glueing} in \cite{glue}) shared by $\spaces$, $\loc$, $\topos$, $\cat$,  and $\pos$.   We will also show that the 1-cotabulators in  $\ring$ are not strong, and so consideration of only strong cotabulators would eliminate this example from consideration.

\smallskip
The cotabulator of  $m\colon X_0\tov X_1$  in $\cat$ (and similarly, $\pos$), also known as the {\it collage}, is the the category $X$ over $\mathbb 2$ whose fibers over $0$ and $1$ are $X_0$ and $X_1$, respectively, and morphism from objects of  $X_0$ to those of $X_1$ are given by $m\colon X_0^{op}\times X_1\to<125>{\rm Sets}$, with the obvious cell $m\to<125>\idv_X$.

\smallskip
Cotabulators in  $\topos$, $\loc$, and   $\spaces$ are constructed using Artin-Wraith glueing (see \cite{TT}, \cite{inclusions}, \cite{glue}). In particular, given $m\colon X_0\tov X_1$ in $\spaces$,  the points of $\Gamma m$ are given by the disjoint union of $X_0$ and $X_1$ with $U$ {\it open} in $\Gamma m$ if and only if $U_0$ is open in $X_0$, $U_1$ is open in $X_1$, and $U_1\subseteq m(U_0)$, where $U_i=U\cap X_i$.

\smallskip
For $\cospan({\cal D})$, the cotabulator of $X_0\to<125>^{c_0}X\toleft<125>^{c_1}X_1$ is given by $X$ with cell $(c_0,\id_X,c_1)$.  If $\cal D$ has pullbacks and pushouts, then the cotabulator of $X_0\toleft<125>^{s_0}X\to<125>^{s_1}X_1$ in $\span({\cal D})$ is the pushout of $s_0$ and $s_1$.  

\smallskip
The situation in $\mathbb M{\rm od}({\cal V})$ is more complicated.  Suppose $M\colon X_0\tov X_1$, i.e., $M$ is an  $(X_0,X_1)$-bimodule.  Then $M$ is an  $X_0\otimes X_1$-module, and so (with appropriate assumptions which apply to $\ring$ and $\quant$), we can consider the symmetric $X_0\otimes X_1$-algebra $SM$, and it is not difficult to show that the inclusion $M\to<125>SM$ defines a cell which gives $SM$ the structure of a 1-cotabulator of $M$.  However, as shown by the following example, the tetrahedron condition need not hold in $\mathbb M{\rm od}({\cal V})$. 

\smallskip
Consider $0\colon \mathbb Z\tov\mathbb Z$ together with $S0=\mathbb Z$ and the unique homomorphism $0\to<125>\mathbb Z$ in $\mathbb M{\rm od}({\rm Ab})=\ring$.  Then, taking $\iota_1(n)=(n,0)$ and  $\iota_2(n)=(0,n)$, the diagram
$$\bfig
\square<300,300>[\mathbb Z`\mathbb Z`\mathbb Z`\mathbb Z;`0`\mathbb Z`]
\square(0,300)<300,300>[\mathbb Z`\mathbb Z`\mathbb Z`\mathbb Z;`\mathbb Z`\mathbb Z\oplus\mathbb Z`]
\place(150,150)[\to<100>^{\scriptscriptstyle0}]
\place(150,450)[\to<100>^{\iota_1}]
\place(600,300)[=]
\square(900,300)<300,300>[\mathbb Z`\mathbb Z`\mathbb Z`\mathbb Z;`0`\mathbb Z`]
\square(900,0)<300,300>[\mathbb Z`\mathbb Z`\mathbb Z`\mathbb Z;`\mathbb Z`\mathbb Z\oplus\mathbb Z`]
\place(1050,450)[\to<100>^{\scriptscriptstyle0}]
\place(1050,150)[\to<100>^{\iota_2}]
\efig$$
defines a commutative tetrahedron which does not factor
$$\bfig
\Dtriangle/->`->`<-/<300,200>[\mathbb Z`\mathbb Z`\mathbb Z;0``]
\Ctriangle(300,0)/<-`->`->/<300,200>[\mathbb Z`\mathbb Z`\mathbb Z;`\mathbb Z\oplus\mathbb Z`]
\place(480,200)[\to<100>^{\varphi}]
\place(120,200)[\to<100>^{\scriptscriptstyle0}]
\place(0,200)[\hbox{\tiny${\bullet\atop{ \ }}$}]
\place(600,200)[\hbox{\tiny${\bullet\atop{ \ }}$}]
\efig$$
for any homomorphism $\varphi\colon \mathbb Z\to<125> \mathbb Z\oplus\mathbb Z$.  Thus, the 1-cotabulators in $\mathbb M{\rm od}({\cal V})$ need not be strong.

\smallskip
To define $\mathbb2$-glueing, suppose $\mathbb D$ has 1-cotabulators and a terminal object $1$, and let $\mathbb 2$ denote the image under $\Gamma$ of the vertical identity morphism on $1$, where $\Gamma\colon\mathbb D_1\to<125>\mathbb D_0$ is left adjoint to $\Delta$ (see Proposition~\ref{tab}).  Then $\Gamma$ induces a functor $\mathbb D_1\to<125>\mathbb D_0/\mathbb 2$, which we also denote by $\Gamma$.  If this functor is an equivalence of categories, then we say {\it $\mathbb D$ has $\mathbb2$-glueing}.  

\smallskip
In $\cat$ (and similarly, $\pos$), $\mathbb 2$ is the  category with two objects and one non-identity morphism.  It is the Sierpinski space $\mathbb 2$ in $\spaces$, the Sierpinski locale ${\cal O}(\mathbb 2)$ in $\loc$, and the Sierpinski topos $S^{\mathbb 2}$ in $\topos$.  That $\cat$ has  $\mathbb2$-glueing is B\'{e}nabou's equivalence cited in \cite{Street}.  For $\spaces$,  $\loc$, and $\topos$, the equivalence follows from the glueing construction (see \cite{TT}, \cite{inclusions}, \cite{glue}).  Note that in each of these cases, $\mathbb 2$ is exponentiable in  $\mathbb D_0$ (see \cite{Hyland,TT}), and so the functor ${\mathbb 2}^*\colon\mathbb D_0 \to<125>\mathbb D_0/\mathbb 2$ has a right adjoint, usually denoted by $\Pi_{\mathbb 2}$.

\begin{prop}\label{glue} If $\mathbb D$ has  $\mathbb2$-glueing and $\mathbb2$ is exponentiable in  $\mathbb D_0$, then $\mathbb D$ has 1-tabulators.
\end{prop}

\noindent{\bf Proof}.  Consider the composite $F\colon\mathbb D_0\to<125>^{\mathbb2^*}\mathbb D_0/\mathbb 2\simeq\mathbb D_1$. Since it is not difficulty to show that $\Gamma$ takes the vertical identity on $X$ to the projection $X\times\mathbb2\to<125>\mathbb2$, it follows that $F=\Delta$, and so $\Delta$ has a right adjoint $\Sigma$, since $\mathbb2^*$ does.  Thus, $\mathbb D$ has 1-tabulators by Proposition~\ref{tab}.
\hfill$\Box$

\medskip
Applying Proposition~\ref{glue}, we see that $\cat$, $\pos$, $\spaces$, $\loc$, and $\topos$ have 1-tabulators (which are can be shown to be strong).  Unraveling the construction of $\Sigma$ given in the proof above, one gets the following descriptions of 1-tabulators in $\cat$ and $\spaces$ which can be shown to be strong.  

\smallskip
Given $m\colon X_0\tov X_1$ in $\cat$ (and similarly $\pos$), the tabulator is the category of elements of $m$, i.e., the objects of $\Sigma m$ are of the form $(x_0,x_1,\alpha)$, where $x_0$ is an object of $X_0$, $x_1$ is an object of $X_1$, and $\alpha\in m(x_0,x_1)$. Morphisms  from $(x_0,x_1,\alpha)$ to $(x'_0,x'_1,\alpha')$ in $\Sigma m$  are pairs $(x_0\to<125>x'_0,x_1\to<125>x'_1)$ of morphisms, compatible with $\alpha$ and $\alpha'$. 

\smallskip
The tabulator of $m\colon X_0\tov X_1$ in $\spaces$ is the set $$\Sigma m=\{(x_0,x_1)\mid \forall U_0\in {\cal O}(X_0), x_1\in m(U_0)\Rightarrow x_0\in U_0\}\subseteq X_0\times X_1$$  
with the subspace topology.  Note that one can directly see that this is the tabulator of $m$ by showing that $(f_0,f_1)\colon Y\to<125>X_0\times X_1$ factors through $\Sigma m$ if and only if $f_1^{-1}m\subseteq f_0^{-1}$.

\smallskip
Although $\span({\cal D})$ and $\cospan({\cal D})$ do not have $\mathbb2$-glueing (since $\Gamma(\idv_1)=1$), and so Proposition~\cite{glue} does not apply, the construction of their tabulators is dual to that of their cotabulators.  Finally,  $\ring$ and $\quant$ do not have 1-tabulators, as they are special cases of Corollary~\ref{tabcor}.

\section{The Adjunction}
Recall from \cite{GP99} that a {\it lax functor} $F\colon \mathbb D\to<125> \mathbb E$ consists of functors $F_i\colon \mathbb D_i\to<125> \mathbb E_i$, for $i=0,1$, compatible with $d_0$ and $d_1$; together with identity and composition comparison cells
$$\rho_X\colon \idv_{FX}\to<125>F(\idv_X)\qquad{\rm and} \qquad   \rho_{m,m'} \colon Fm'\dotv Fm\to<125> F(m'\dotv m)$$  for every object $X$ and every vertical composite $m'\dotv m$ of $\mathbb D$, respectively;  satisfying naturality and coherence conditions.  If $\rho_X$ is an isomorphism, for all $X$, we say that $F$ is a {\it normal lax functor}.  An {\it oplax functor}  is defined dually with comparison cells in the opposite direction.

\smallskip
An {\it oplax/lax adjunction} consists of an oplax functor $F\colon \mathbb D\to<125> \mathbb E$ and a lax double functor $G\colon \mathbb E\to<120> \mathbb D$ together with double cells 
$$\bfig
\square<400,350>[X_0`GFX_0`X_1`GFX_1;\eta_{X_0}`m`GFm`\eta_{X_1}]
\place(0,200)[\hbox{\tiny$\bullet$}]
\place(400,200)[\hbox{\tiny$\bullet$}]
\place(200,175)[\eta_m]
\square(1200,0)<400,350>[FGY_0`Y_0`FGY_1`Y_1;\varepsilon_{Y_0}`FGn`n`\varepsilon_{Y_1}]
\place(1200,200)[\hbox{\tiny$\bullet$}]
\place(1600,200)[\hbox{\tiny$\bullet$}]
\place(1400,175)[\varepsilon_n]
\efig$$
satisfying naturality  and coherence conditions, as well as the usual adjunction identities (see  \cite{GP04}).

\begin{ex}
Suppose $\mathbb  D$ is a double category with 1-cotabulators and $\mathbb D_0$ has pushouts.  Then, by Proposition~\ref{tab},  the functor $\Delta\colon \mathbb D_0\to<125> \mathbb D_1$ has a left adjoint (denoted by $\Gamma$), and so there is an oplax functor $F\colon \mathbb D\to<125>{\cospan}(\mathbb D_0)$ which is the identity on objects and horizontal morphisms, and defined on vertical morphisms and cells by  
$$\bfig
\square<350,350>[X_0`X'_0`X_1`X'_1;f_0`m`m'`f_1]
\place(0,200)[\hbox{\tiny$\bullet$}]
\place(350,200)[\hbox{\tiny$\bullet$}]
\place(175,175)[\varphi]
\place(800,200)[\longmapsto]
\Dtriangle(1200,0)|rba|/`->`<-/<400,200>[X_0`\Gamma m`X_1;`i_0`i_1]
\Dtriangle(1700,0)/`->`<-/<400,200>[X'_0`\Gamma  m'`X'_1;`i'_0`i'_1]
\morphism(1200,400)<500,0>[X_0`X'_0;f_0]
\morphism(1200,0)|b|<500,0>[X_1`X'_1;f_1]
\morphism(1600,200)<500,0>[\Gamma m`\Gamma m';]
\place(1800,260)[\scriptstyle f]
\efig$$
where $f$ is induced by the universal property of the 1-cotabulator.  The comparison cells $F(\idv_X)\to<150>\idv_{FX}$  and $F(m'\dotv m)\to<150>Fm'\dotv Fm$ also arise via the universal property, with the latter given by the horizontal morphism $\Gamma(m'\dotv m)\to<150>P$  corresponding to the diagram
$$\bfig
\Dtriangle(400,0)|rba|/`->`<-/<400,200>[X_1`\Gamma m'`X_2;``]
\Dtriangle(400,400)|rba|/`->`<-/<400,200>[X_0`\Gamma m`X_1;``]
\morphism(0,800)<400,0>[X_0`X_0;\id_{X_0}]
\morphism(0,0)|b|<400,0>[X_2`X_2;\id_{X_2}]
\morphism(0,800)<0,-800>[X_0`X_2;m'\dotv m]
\morphism(400,800)<0,-400>[X_0`X_1;m]
\morphism(400,400)<0,-400>[X_1`X_2;m']
\place(0,400)[\hbox{\tiny$\bullet$}]
\place(400,200)[\hbox{\tiny$\bullet$}]
\place(400,600)[\hbox{\tiny$\bullet$}]
\place(530,600)[\iota_m]
\place(530,200)[\iota_{m'}]
\morphism(800,600)<400,-200>[\Gamma m`P;]
\morphism(800,200)<400,200>[\Gamma m'`P;]
\efig$$
where $P$ is a pushout and the large rectangle commutes.
\end{ex}

Dually, we get:

\begin{ex} Suppose $\mathbb D$ is a double category with 1-tabulators and $\mathbb D_0$ has pullbacks.  Then there is a 
lax functor $F\colon \mathbb D\to<125>{\span}(\mathbb D_0)$ which is the identity on objects and horizontal morphisms and takes $m\colon X_0\tov X_1$ to the span $X_0\to/<-/<125> \Sigma m\to<125> X_1$, where $\Sigma$ is the right adjoint to $\Delta$.
\end{ex}

\begin{prop}\label{framed} 
Suppose $\mathbb D$ is a double category and $\mathbb D_0$ has pushouts.  Then $\mathbb D$ has companions and conjoints if and only the identity functor on $\mathbb D_0$ extends to a normal lax functor $G\colon{\cospan}(\mathbb D_0)\to<125> \mathbb D$.
\end{prop}

\noindent{\bf Proof}.  Suppose $\mathbb D$ has companions and conjoints.  Then it is not difficult to show that there is  a normal lax functor $G\colon {\cospan}(\mathbb D_0)\to<125> \mathbb D$ which is the identity on objects and horizontal morphisms, and is defined on cells by
$$\bfig
\Dtriangle|lll|/`->`<-/<350,200>[Y_0`Y`Y_1;`c_0`c_1]
\morphism(350,200)<500,0>[Y`Y';g]
\morphism(0,400)|a|<500,0>[Y_0`Y'_0;g_0]
\morphism(0,0)|b|<500,0>[Y_1`Y'_1;g_1]
\Dtriangle(500,0)/`->`<-/<350,200>[Y'_0`Y'`Y'_1;`c'_0`c'_1]
\place(1150,200)[\longmapsto]
\square(1500,200)|alra|<350,300>[Y_0`Y'_0`Y`Y';g_0`{c_0}_*`{c'_0}_*`g] 
\square(1500,-100)<350,300>[Y`Y'`Y_1`Y'_1;`c_1^*`{c'_1}^{\!\!*}`g_1] 
\place(1500,400)[\hbox{\tiny$\bullet$}]
\place(1850,400)[\hbox{\tiny$\bullet$}]
\place(1500,100)[\hbox{\tiny$\bullet$}]
\place(1850,100)[\hbox{\tiny$\bullet$}]
\place(1690,400)[\scriptstyle{\psi_0}]
\place(1690,70)[\scriptstyle{\psi_1}]
\efig$$
where $\psi_0$ and $\psi_1$ arise from the commutativity of the squares in the cospan cell, and the definitions of companion and conjoint.

\smallskip
Conversely, suppose there is a normal lax functor $G\colon{\cospan}(\mathbb D_0)\to<125> \mathbb D$ which is the identity on objects and horizontal morphisms. Then the companion and conjoint of $f\colon X\to<125>Y$ are defined as follows.  Consider $$f_*=G(X\to<125>^fY\to/<-/<125>^{\id_Y} Y)\colon X\tov Y \hskip .5in f^*=G(Y\to<125>^{\id_Y} Y\to/<-/<125>^fX)\colon Y\tov X$$

\smallskip
Applying $G$ to the cospan diagrams
$$\bfig
\Dtriangle|lll|/`->`<-/<350,200>[X`X`X;`\id_X`\id_X]
\morphism(0,400)|a|<400,0>[X`X;\id_X]
\morphism(350,200)<400,0>[X`Y;\!\!\!\!f]
\morphism(0,0)|b|<400,0>[X`Y;f]
\Dtriangle(400,0)/`->`<-/<350,200>[X`Y`Y;`f`\id_Y]
\Dtriangle(1200,0)|lll|/`->`<-/<350,200>[X`Y`Y;`f`\id_Y]
\morphism(1200,400)|a|<400,0>[X`Y;f]
\morphism(1550,200)<400,0>[Y`Y;\!\!\!\!\id_Y]
\morphism(1200,0)|b|<400,0>[Y`Y;\id_Y]
\Dtriangle(1600,0)/`->`<-/<350,200>[Y`Y`Y;`\id_Y`\id_Y]
\efig$$
and composing with $\rho_X$ and $\rho_Y^{-1}$, we get cells
$$\bfig
\square<350,350>[X`X`X`X;\id_X`\idv_X``\id_X]
\place(0,200)[\hbox{\tiny$\bullet$}]
\place(350,200)[\hbox{\tiny$\bullet$}]
\place(120,150)[\scriptstyle{\rho_X}]
\square(350,0)<375,350>[X`X`X`Y;\id_X`{\id_X}_{\!*}`f_*`f]
\place(725,200)[\hbox{\tiny$\bullet$}]
\place(535,150)[\scriptscriptstyle{G(\id_X,f,f)}]
\square(1300,0)<400,350>[X`Y`Y`Y;f`f_*`{\id_Y}_{\!*}`\id_Y]
\place(1300,200)[\hbox{\tiny$\bullet$}]
\place(1700,200)[\hbox{\tiny$\bullet$}]
\place(1500,150)[\scriptscriptstyle{G(f,\id_y,\id_Y)}]
\square(1700,0)<350,350>[Y`Y`Y`Y;\id_Y``\idv_Y`\id_Y]
\place(2050,200)[\hbox{\tiny$\bullet$}]
\place(1940,150)[\scriptstyle{\rho_Y^{-1}}]
\efig$$
which serve as $\eta$ and $\varepsilon$, respectively, making $f_*$ the companion of $f$.  Note that the horizontal and vertical identities for $\eta$ and $\varepsilon$ follow from the normality and coherence axioms of $G$, respectively.

\smallskip
Similarly, the cells $\alpha$ and $\beta$ for $f^*$ arise from the cospan diagrams 
$$\bfig
\Dtriangle|lll|/`->`<-/<350,200>[X`X`X;`\id_X`\id_X]
\morphism(0,400)|a|<400,0>[X`Y;f]
\morphism(350,200)<400,0>[X`Y;\!\!\!\!f]
\morphism(0,0)|b|<400,0>[X`X;\id_X]
\Dtriangle(400,0)/`->`<-/<350,200>[Y`Y`X;`\id_Y`f]
\Dtriangle(1200,0)|lll|/`->`<-/<350,200>[Y`Y`X;`\id_Y`f]
\morphism(1200,400)|a|<400,0>[Y`Y;\id_Y]
\morphism(1550,200)<400,0>[Y`Y;\!\!\!\!\id_Y]
\morphism(1200,0)|b|<400,0>[X`Y;f]
\Dtriangle(1600,0)/`->`<-/<350,200>[Y`Y`Y;`\id_Y`\id_Y]
\efig$$
and it follows that $\mathbb D$ has companions and conjoints.
\hfill$\Box$

\begin{cor}
Suppose $\mathbb D$ is a double category and $\mathbb D_0$ has pullbacks.  Then $\mathbb D$ has companions and conjoints if and only the identity functor on $\mathbb D_0$ extends to an opnormal oplax functor ${\span}(\mathbb D_0)\to<125> \mathbb D$.
\end{cor}

\noindent{\bf Proof}.  Apply Proposition~\ref{framed} to $\mathbb D^{op}$.
\hfill$\Box$

\begin{thm}\label{adj} The following are equivalent for a double category $\mathbb D$ such that $\mathbb D_0$ has pushouts:
\begin{description}
\item{(a)} there is an oplax/lax adjunction $\mathbb D\two/->`<-/<150>^F_G\cospan(\mathbb D_0)$ such that $G$ is normal and restricts to the identity on $\mathbb D_0$;
\item{(b)} $\mathbb D$ has companions, conjoints, and 1-cotabulators;
\item{(c)} $\mathbb D$ has companions and conjoints, and the induced normal lax functor $G\colon \cospan(\mathbb D_0) \to<125>\mathbb D$ has an oplax left adjoint;
\item{(d)} $\mathbb D$ has companions, conjoints, and 1-cotabulators, and the induced oplax functor $F\colon\mathbb D \to<125>\cospan(\mathbb D_0)$ is left adjoint to the induced normal lax functor $G\colon \cospan(\mathbb D_0) \to<125>\mathbb D$;
\item{(e)} $\mathbb D$ has 1-cotabulators and the induced oplax functor $F\colon\mathbb D \to<125>\cospan(\mathbb D_0)$ has a normal lax right adjoint. 
\end{description}
\end{thm}

\noindent{\bf Proof}. (a)$\Rightarrow$(b) Given $F\dashv G$ such that  $G$ is normal and restricts to the identity on $\mathbb D_0$, we know $\mathbb D$ has companions and conjoints by Proposition~\ref{framed}.  To see that $\mathbb D$ has 1-cotabulators, by Proposition~\ref{tab}, it suffices to show that $\Delta\colon \mathbb D_0\to<125> \mathbb D_1$ has a left adjoint.  Since $\Delta$ factors as $\mathbb D_0\to<125>^\Delta{\rm Cospan}(\mathbb D_0)\to<125>^{G_1} \mathbb D_1$, by normality of $G$, and both these functors have left adjoints,  the desired result follows.

\smallskip
\noindent (b)$\Rightarrow$(c) Suppose $\mathbb D$ has companions, conjoints, and 1-cotabulators, and consider the induced functor $F\colon\mathbb D \to<125>\cospan(\mathbb D_0)$.  Given $m\colon X_0\tov X_1$ and $Y_0\to<125>^{c_0}Y\toleft<125>^{c_1}Y_1$, applying the definition of 1-cotabulator, we know that every cell in $\cospan(\mathbb D_0)$ of the form
$$\bfig
\Dtriangle|lll|/`->`<-/<350,200>[X_0`\Gamma m`X_1;``]
\morphism(0,400)|a|<400,0>[X_0`Y_0;f_0]
\morphism(350,200)<400,0>[\Gamma m`Y;]
\morphism(0,0)|b|<400,0>[X_1`Y_1;f_1]
\Dtriangle(400,0)/`->`<-/<350,200>[Y_0`Y`Y_1;`c_0`c_1]
\efig$$
corresponds to a unique cell 
$$\bfig
\square<350,350>[X_0`Y`X_1`Y;c_0f_c`m`\id_Y`c_1f_1]
\place(175,175)[\varphi]
\place(0,175)[\hbox{\tiny${\bullet\atop{ \ }}$}]
\place(350,175)[\hbox{\tiny${\bullet\atop{ \ }}$}]
\efig$$
and hence, a unique cell
$$\bfig
\square(0,0)/->`->``->/<350,600>[X_0`Y_0`X_1`Y_1;f_0`m``f_1]
\morphism(350,600)|r|<0,-300>[Y_0`Y;{c_0}_*]
\morphism(350,300)|r|<0,-300>[Y`Y_1;c_1^*]
\place(0,300)[\hbox{\tiny${\bullet\atop{ \ }}$}]
\place(350,150)[\hbox{\tiny${\bullet\atop{ \ }}$}]
\place(350,450)[\hbox{\tiny${\bullet\atop{ \ }}$}]
\place(175,300)[\psi]
\efig$$
by vertical flipping, and it follows that  $F\dashv G$.

\smallskip
\noindent (c)$\Rightarrow$(d) Suppose $\mathbb D$ has companions and conjoints, and the induced normal lax functor $G$ has an oplax left adjoint $F$.  As in the proof of (a)$\Rightarrow$(b), we know that $\mathbb D$ has 1-cotabulators, and so, it suffices to show that $F$ is the induced functor.  We know $F$ takes $m\colon X_0\tov X_1$ to a cospan of the form $X_0\to<125>^{c_0}\Gamma m\toleft<125>^{c_1}X_1$, since $F$ is the identity on objects and the left adjoint of $\Delta\colon \mathbb D_0\to<125>{\rm Cospan}(\mathbb D_0)$ takes $X_0\to<125>^{c_0}\Gamma m\toleft<125>^{c_1}X_1$ to $X$, and so the desired result easily follows.  

\smallskip
\noindent (d)$\Rightarrow$(e) is clear.

\smallskip
\noindent (e)$\Rightarrow$(a) Suppose $\mathbb D$ has 1-cotabulators and the induced oplax functor $F$ has a normal lax right adjoint $G$.  Since $F$ restricts to the identity on $\mathbb D_0$, then so does $G$, and the proof is complete.
\hfill$\Box$

\medskip
Note that this proof shows that if there is an oplax/lax adjunction $$\mathbb D\two/->`<-/<150>^F_G\cospan(\mathbb D_0)$$ such that $G$ is normal and restricts to the identity on $\mathbb D_0$, then $F$ is the induced by 1-cotabulators and $G$  by companions and conjoints.  Since the double categories in Examples~2.1--2.8 all have companions, conjoints, and 1-cotabulators, it follows that they each admits a unique (up to equivalence) oplax/lax adjunction of this form and it is induced in this manner.

\smallskip
Applying Theorem~\ref{adj} to $\mathbb D^{op}$, we get:

\begin{cor}\label{adjop}  The following are equivalent for a double category $\mathbb D$ such that $\mathbb D_0$ has pullbacks:
\begin{description}
\item{(a)} there is an oplax/lax adjunction $\span(\mathbb D_0)\two/->`<-/<150>^G_F\mathbb D$ such that $G$ is opnormal and restricts to the identity on $\mathbb D_0$;
\item{(b)} $\mathbb D$ has companions, conjoints, and 1-tabulators;
\item{(c)} $\mathbb D$ has companions and conjoints, and the induced opnormal oplax functor $G\colon \span(\mathbb D_0) \to<125>\mathbb D$ has an lax right adjoint;
\item{(d)} $\mathbb D$ has companions, conjoints, and 1-tabulators, and the induced lax functor $F\colon\mathbb D \to<125>\span(\mathbb D_0)$ is right adjoint to the induced opnormal oplax functor $G\colon \span(\mathbb D_0) \to<125>\mathbb D$;
\item{(e)} $\mathbb D$ has 1-tabulators and the induced lax functor $F\colon\mathbb D \to<125>\span(\mathbb D_0)$ has a opnormal oplax left adjoint. 
\end{description}
\end{cor}

\medskip
As in the cospan case, if there is an oplax/lax adjunction $$\span(\mathbb D_0)\two/->`<-/<150>^G_F\mathbb D$$ such that $G$ is opnormal and restricts to the identity on $\mathbb D_0$, then $F$ is the induced by 1-tabulators and $G$ is by companions and conjoints.  Since the double categories in Examples~2.1--2.7 (i.e., all by $\mathbb M{\rm od}({\cal V})$ ) have companions, conjoints, and 1-tabulators, it follows they each admits is a unique (up to equivalence) oplax/lax adjunction of this form, and it is induced by companions, conjoints, and 1-tabulators.

\end{document}